\documentclass[12pt,a4paper]{amsart}
\usepackage{amssymb, amsthm, amsbsy}

\newtheorem{theor}{Theorem}

\newtheorem{state}[theor]{Proposition}

\theoremstyle{definition}
\newtheorem{define}{Definition}
\newtheorem{cor}[theor]{Corollary}

\theoremstyle{remark}
\newtheorem{rem}{Remark}
\newtheorem{example}{Example}

\newcommand{\Id}{{\mathrm d}}
\newcommand{\ID}{{\mathrm D}}

\newcommand{\pinner}{\mathbin{\mathchoice
   {\hbox{\vrule width0.6em depth0pt height0.4pt
   \vrule width0.4pt depth0pt height0.8ex}}
   {\hbox{\vrule width0.6em depth0pt height0.4pt
   \vrule width0.4pt depth0pt height0.8ex}}
   {\hbox{\kern0.14em
   \vrule width0.48em depth0pt height0.4pt
   \vrule width0.4pt depth0pt height0.6ex\kern0.14em}}
   {\hbox{\kern0.1em
   \vrule width0.39em depth0pt height0.4pt
   \vrule width0.4pt depth0pt height0.5ex\kern0.1em}}}}
\newcommand{\inner}{\pinner\,}

\DeclareFontFamily{OML}{cyr}{}
\DeclareFontShape{OML}{cyr}{m}{n}{
   <5> <6> <7> <8> <9> gen * wncyr
   <10> <10.95> <12> <14.4> <17.28> <20.74> <24.88> wncyr10
  }{}
\DeclareSymbolFont{rusletters}{OML}{cyr}{m}{n}
\DeclareSymbolFontAlphabet{\rusmath}{rusletters}
\DeclareMathSymbol\re{\rusmath}{rusletters}{"03}

\newcommand{\cEv}{\re}

\newcommand{\BBC}{{\mathbb{C}}}

\newcommand{\BBR}{{\mathbb{R}}}
\newcommand{\BBN}{{\mathbb{N}}}
\newcommand{\BBZ}{{\mathbb{Z}}}

\newcommand{\const}{\mathop{\rm const}\nolimits}

\newcommand{\bx}{{\boldsymbol{x}}}

\newcommand{\cA}{\mathcal{A}}

\newcommand{\cC}{\mathcal{C}}

\newcommand{\cF}{\mathcal{F}}

\newcommand{\cO}{\mathcal{O}}

\newcommand{\dd}{\partial}

\newcommand{\bun}{\mathbf{1}}

\newcommand{\nel}{{\text{\textup{in}}}}
\newcommand{\out}{{\text{\textup{out}}}}

\DeclareMathOperator{\CDiff}{\mathcal{C}Diff}

\DeclareMathOperator{\Lie}{Lie}
\DeclareMathOperator{\Hom}{Hom}
\DeclareMathOperator{\Diff}{Diff}

\DeclareMathOperator{\charact}{char}

\newcommand{\sch}[2]{{[\![{#1},{#2}]\!]}^{\mathrm{S}}}
\newcommand{\rn}[2]{{[\![{#1},{#2}]\!]}^{\mathrm{RN}}}
\newcommand{\W}[3]{\mathop{{W}}\left(#1,#2,#3\right)}

\newcommand{\by}[1]{\textrm{{#1}}}
\newcommand{\jour}[1]{\textit{{#1}}}
\newcommand{\vol}[1]{\textbf{{#1}}}
\newcommand{\book}[1]{\textit{{#1}}}

\title[On the associative homotopy Lie algebras]%
{On the associative homotopy Lie algebras and the Wronskians}

\date{October 28, 2003}

\author{Arthemy V. Kiselev}
\thanks{Submitted to: \textit{Fundamental'naya i Prikladnaya Matematika}
(2004), English translation: \textit{Journal of Mathematical Sciences}}
\thanks{Partially supported by the INTAS grant YS~2001/2-33}

\address{Ivanovo State Power University, Dept.\ of Higher Mathematics,
Rabfakovskaya str.\ 34, Ivanovo, 153003 Russia.}

\curraddr{Department of Mathematics, Brock University, 500 Glenridge Ave.,
St.~Catharines, Ontario, Canada L2S~3A1.}

\email{Arthemy.Kiselev@BrockU.CA}



\begin{document}
\rightline{ISPUmath-02/2004}

\begin{abstract}
Representations of the Schlessinger\/--\/Stasheff's
associative homotopy Lie algebras in the spaces of
higher\/--\/order differential operators are
analyzed. The $W$-transformations of chiral
embeddings, related with the Toda equations,
of complex curves into the K\"ahler
manifolds are shown to be endowed with the homotopy Lie algebra
structures.
Extensions of the Wronskian determinants that preserve
the properties of the Schlessinger\/--\/Stasheff algebras
are constructed for the case of $n\geq1$
independent variables.
\end{abstract}

\maketitle
\tableofcontents

\subsection*{Introduction}
%
Recently, the attention of the mathematical physics community has been
drawn to the $N$-ary Lie algebra structures, i.e., the $N$-linear
skew\/--\/symmetric brackets that satisfy an analog of the Jacobi
identity, and to the $N$-field dynamics problems, e.g., $N$-ary objects
similar to the Poisson manifolds.
There are several concepts (\cite{DzhFr, Filippov, HanlonWachs,
Stasheff, Sahoo, SS}) how the Jacobi identity should be generalized,
each of them having its own interpretation (e.g., \cite{Nambu}) and
applications (\cite{DzhFAP}). Recently, a unifying approach was
proposed in \cite{Uvas}, treating the known cases as the special
ones within the $3$-parameter family of the identities.

Up to our present knowledge, the papers \cite{Filippov, SS}, issued in
1985, were the first works on the topic. In \cite{Filippov}, V. T.
Filippov considered an $N$-linear skew\/--\/symmetric bracket $\nabla$,
defined on a vector space $\cA$, that satisfied an analog of the Jacobi
identity:
\begin{equation}\tag{I.1}\label{J1}
\nabla\left(a_1,\ldots,a_{N-1},\nabla(b_1,\ldots,b_N)\right)=
\sum_{i=1}^N
\nabla\left(b_1,\ldots,b_{i-1},
\nabla(a_1,\ldots,a_{N-1},b_i),b_{i+1},\ldots,b_N\right),
\end{equation}
where $a_i$, $b_j\in\cA$; then, $\cA$ is called an $N$-Lie algebra.
The motivation of \eqref{J1} to appear is quite understandable: the
adjoint representation $a\mapsto\nabla(a_1$, $\ldots$, $a_{N-1}$, $a)$
is a derivation for any $a_i\in\cA$. The standard constructions of the
Lie algebra theory for the $N$-Lie algebras were introduced in
\cite{Kasymov}.

The Nambu mechanics is an example of the $N$-Poisson
dynamics assigned to identity \eqref{J1}; in \cite{Nambu}, the standard
binary Poisson bracket was replaced by the ternary ($N=3$) one:
\[
\nabla(f_1,\dots,f_N)=\det\left\|\frac{\dd f_i}{\dd x^j}\right\|,
\]
where $\cA=C^\infty(\BBR^N)$ and the r.h.s. contains the Jacobi
determinant; nevertheless, the fact that Eq.\ \eqref{J1} holds for this
$\nabla$ was not noticed until \cite{Sahoo}. Then, L. Takhtajan
(\cite{Takhtajan}) developed the concept of the Nambu\/--\/Poisson
manifolds for $N\geq2$.

\medskip
The second natural generalization of the ordinary Jacobi identity is
\begin{equation}\tag{I.2}\label{J2}
\sum\limits_{\sigma\in S^N_{2N-1} } (-1)^\sigma\Delta(\Delta(
a_{\sigma(1)},\ldots,
   a_{\sigma(N)}), a_{\sigma({N+1})},\ldots, a_{\sigma(2N-1)}) = 0,
\end{equation}
where $a_i\in\cA$ and $S_{2N-1}^N=\{\sigma\in
S_{2N-1}\mid\sigma(1)<\dots<\sigma(N)$,
$\sigma(N+1)<\dots<\sigma(2N-1)\}$ is the set of
the unshuffle permutations.
These brackets $\Delta$ are named the homotopy $N$-Lie algebra
structures (\cite{SS}) and are closely related with the
SH\/--\/algebras (\cite{Barnich, Stasheff}). Also, these algebras and
their Koszul cohomologies were studied in \cite{HanlonWachs}; their
Hochschild cohomologies were considered in \cite{MichorVin}.

The $N$-Poisson manifolds associated with Jacobi's identity \eqref{J2}
were introduced in \cite{Perelomov}: an $N$-vector field $V$ is an
$N$-Poisson structure if the equation $\sch{V}{V}=0$ for the Schouten
bracket holds (see page \pageref{RNBracket} for definitions).

The properties of the Filippov's $N$-Lie and the
Schlessinger\/--\/Stasheff's homotopy $N$-Lie algebras are quite
different; really, they appear in the $3$-parametric scheme
(\cite{Uvas}) as the opposite cases: $(N$, $N-1$, $0)$ and $(N$, $0$,
$0)$, respectively (see Definition \ref{Defnkr} below).
Further discussion on the topic is found in \cite{Uvas}.

\bigskip
In the sequel, we analyze the properties of \emph{associative}
homotopy Lie algebras and their representations in higher order
differential operators. Namely, we relate the corresponding structures
with the Wronskian determinants, point out $N$-ary analogues of
vector fields on smooth manifolds, and construct a definition of the
Wronskian for $n\geq1$ independent variables which is correlated with
the structures defined by Eq.\ \eqref{J2}.  Also, we prove that
the $W$-transformations of the $W$-surfaces
(\cite{GervaisLetter}--\cite{GervaisENS}) in the K\"ahler manifolds are
endowed with the homotopy Lie algebra structures.
We use the jet bundles language \cite{ClassSymEng};
our approach is aimed to
contribute the study of related aspects in the
cohomological algebra and the field theory.
The exposition patterns upon \cite[Chapter~V]{Thesis}.

The paper is organized as follows.

In Section \ref{SecAlgConcept},
we introduce the main algebraic concept of the
homotopy $N$-Lie algebras. We fix notation and define the
Richardson\/--\/Nijenhuis bracket $\rn{\cdot}{\cdot}$,
the homotopy $N$-Jacobi identities $\rn{\Delta}{\Delta}=0$ for
$N$-linear skew\/--\/symmetric operators $\Delta$, and the Hochschild
and Koszul cohomologies.
Next, we illustrate the definitions by two
finite\/--\/dimensional homotopy Lie algebras
which are analogues the Lie algebra $\mathfrak{sl}_2(\Bbbk)$.

In Section \ref{SecAssocAlg}, we consider representations of the
homotopy Lie algebra structures in the higher order differential
operators over $\Bbbk$. Analyzing the properties of the corresponding
$N$-linear skew\/--\/symmetric bracket,
we point out higher\/--\/order
generalizations of vector fields and thus explain the property
$\dd^{\vec\imath}[\dd^{\vec\imath}]=0$ of the Wronskian determinants
$\dd^{\vec\imath}=\bun\wedge\dd\wedge\dots\wedge\dd^{N-1}$ to provide
the homotopy Lie algebra structures for even $N$s;
also, we calculate the structural constants of these algebras in terms
of the Vandermonde determinants.
Next, we relate the multilinear homotopy structures with the Toda
equations (\cite{LeznovSavelievBV}).
The latter are known to be the compatibility
condition in the $W$-geometry (\cite{GervaisLetter}--\cite{GervaisENS})
of chiral embeddings of
complex curves into the K\"ahler manifolds,
while higher order differential operators, endowed with the
homotopy Lie algebra structures, are admissible deformations of
these embeddings.

In Section \ref{SecnDimBase}, we construct analogues
$D^{\vec\sigma}=D_{\sigma_1}\wedge\dots\wedge D_{\sigma_N}$ of the
Wronskian determinants $\dd^{\vec\imath}$ for the case of $n\geq1$
independent variables $x^1$, $\ldots$, $x^n$, such that the
$\binom{n+k}{n}$-Jacobi identities $D^{\vec\sigma}[D^{\vec\sigma}]=0$
are preserved for any integers $n$, $k\geq1$.

\section{Preliminaries: the algebraic concept}\label{SecAlgConcept}
\subsection{Basic definitions and facts}
First let us introduce some notation.
Let $\cA$ be an 
algebra over the field $\Bbbk$ such that $\charact
\Bbbk=0$ and let $\dd$ be a derivation of $\cA$.
As an illustrative example, one can think $\cA$ to be the algebra of
smooth functions $f\colon M\to\BBR$ on a smooth real manifold~$M$.

Let $S^k_m\subset S_m$ be the \emph{unshuffle permutations}
such that $\sigma(1)<\sigma(2)<\cdots
<\sigma(k)$ and $\sigma(k+1)<\sigma(k+2)<\cdots<\sigma(m)$
for any $\sigma\in S^k_m$.
Let $\Delta\in\Hom_\Bbbk(\bigwedge^k\cA,\cA)$
be a homomorphism and take arbitrary $a_j\in\cA$,
$1\leq j\leq k$; suppose $1\leq l\leq k$.
The \emph{inner product}
$\Delta_{a_1,\ldots,a_m}\in\Hom(\bigwedge^{k-m}\cA,\cA)$ is
the operator defined by the rule
$$
\Delta_{a_1,\ldots,a_m}(a_{m+1},\ldots,a_k)\stackrel{\mathrm{def}}{=}
\Delta(a_1,\ldots,a_k).
$$


\begin{define}[\textup{\cite{GelfandFuks68, Uvas}}]\label{exterior}
Let $\Delta\in\Hom(\bigwedge^k\cA,\cA)$ and
$\nabla\in\Hom(\bigwedge^l\cA,\cA)$ be two operators; by definition,
the \emph{exterior multiplication} $\wedge$ in
$\Hom(\bigwedge^*\cA,\cA)$ is
$$
(\Delta\wedge\nabla)(a_1,\ldots,a_{k+l})=\sum\limits_{\sigma\in
S^k_{k+l}}
(-1)^\sigma\,\Delta(a_{\sigma(1)},\ldots,a_{\sigma(k)})\cdot
\nabla(a_{\sigma(k+1)},\ldots,a_{\sigma(k+l)})
$$
for any $a_1,\ldots,a_{k+l}\in\cA$.
\end{define}

\begin{example}
The exterior multiplication $\wedge$ on the space of higher\/--\/order
differential operators that act on the algebra $\cA=C^\infty(M)$ of
smooth functions on $M$ defines the Wronskian determinants
$W^{0,1,\ldots,N+1}=\dd^0\wedge\ldots\wedge\dd^{N+1}$.
In this paper, we also
consider generalized Wronskians $W^{\vec\imath} =
\dd^{i_1}\wedge\ldots\wedge\dd^{i_N}\in\Hom(\bigwedge\nolimits
^N\cA,\cA)$, where $\dd$ is a derivation of~$\cA$. Let a
multiindex $\vec\imath\in\BBZ^N_{+}$ be such that $0\leq
i_1<\ldots<i_N$.
By $\Hom_t(\bigwedge\nolimits
^N\cA,\cA)$ we denote the linear span of the generalized
Wronskians $W^{\vec\imath}$ such that $|\vec\imath|\equiv\sum_ji_j=t$.
By definition, put $|W^{\vec\imath}|=|\vec\imath|$.
In Section \ref{nDimWronsk}, we construct generalizations of the
Wronskian determinants for analytic functions $\Bbbk[[x^1,\ldots,x^n]]$
by using Definition~\ref{exterior}.
\end{example}

Next, let $\Delta\in\Hom(\bigwedge\nolimits^k\cA,\cA)$ and
$\nabla\in\Hom(\bigwedge\nolimits^l\cA,\cA)$. By
$\Delta[\nabla]\in\Hom(\bigwedge^{k+l-1}\cA$, $\cA)$
we denote the \emph{action}
$\Delta[{\cdot}]\colon\Hom(\bigwedge^N\cA,\cA)$ $\to$
$\Hom(\bigwedge^{N+k-1}\cA,\cA)$ of $\Delta$ on $\nabla$:
\begin{equation}\label{Action}
\Delta[\nabla](a_1,\ldots,a_{k+l-1})\stackrel{\mathrm{def}}{{}={}}
\sum\limits_{\sigma\in
S^l_{k+l-1} } (-1)^\sigma\,\Delta(\nabla(a_{\sigma(1)},\ldots,
a_{\sigma(l)}),a_{\sigma(l+1)},\ldots,a_{\sigma(k+l-1)}),
\end{equation}
where $a_j\in\cA$.
By $\rn{\Delta}{\nabla}\in\Hom(\bigwedge^{k+l-1}\cA,\cA)$ we denote the
\emph{Richardson\/--\/Nijenhuis bracket} of $\Delta$ and $\nabla$:
\begin{equation}\label{RNBracket}
\rn{\Delta}{\nabla}\stackrel{\mathrm{def}}{=} \Delta[\nabla] -
(-1)^{(k-1)(l-1)}\nabla[\Delta].
\end{equation}

\begin{define}[\cite{Uvas}]\label{Defnkr}
Choose integers $N$, $k$, and $r$ such that $0\leq r\leq k<N$, and let
$a_1,\ldots,a_r,b_1,\ldots,b_k\in\cA$. The skew-symmetric map
$\Delta\in\Hom(\bigwedge^N\cA,\cA)$ is said to determine the Lie
algebra structure of the type $(N,k,r)$ on the $\Bbbk$-vector space
$\cA$ if $\Delta$ satisfies the $(N,k,r)$-\emph{Jacobi identity}
\begin{equation}\label{nkrJacobi}
\rn{\Delta_{a_1,\ldots,a_r}}{\Delta_{b_1,\ldots,b_k}}=0
\end{equation}
for any $\vec a$ and $\vec b$. By $\Lie^{(N,k,r)}(\cA)$ we denote
the set of all type $(N,k,r)$ structures
$\Delta\in\Hom(\bigwedge^N\cA,\cA)$ on $\cA$.
\end{define}

The structure $\Delta\in\Hom(\bigwedge^N\cA,\cA)$ is called a
\emph{multi-derivation} if the Leibnitz rule
$$
\Delta(ab,a_2,\ldots,a_N)=a\,\Delta(b,a_2,\ldots,a_N) +
\Delta(a,a_2,\ldots,a_N)\,b
$$
is valid for any $a$, $b$, $a_j\in\cA$.

\begin{example}[\textup{Filippov's $N$-Lie algebras}]
Consider the family $\Lie^{(N,N-1,0)}$ for integer $N\geq2$. The
$N$-Jacobi identity is \eqref{J1},
meaning that the adjoint representation for these algebras is a
derivation. This case is a natural generalisation of the Poisson theory
(\cite{Takhtajan}).
\end{example}

\begin{rem}[\textup{\cite{Uvas}}]\label{EquivHeredity}
We have
\begin{equation}\label{HeredityIso}
\Lie^{(N,0,0)}(\cA)=\Lie^{(N,1,0)}(\cA)
\end{equation}
for any even $N$;
this is a typical instance of the heredity structures.
Indeed, the following two conditions are equivalent:
$$
\rn{\Delta}{\Delta}=0\,\Leftrightarrow\,
{\rn{\Delta}{\Delta}}_a=-2\,\rn{\Delta}{\Delta_a}=0\qquad\forall
a\in\cA,
$$
owing to Corollary 1.1 in \cite{Uvas}:
$$
{\rn{\Delta}{\Delta}}_a=(-1)^{N-1}\rn{\Delta}{\Delta_a} +
\rn{\Delta_a}{\Delta}.
$$
Finally, $\rn{\Delta}{\Delta_a}=0$ for any $a\in\cA$.
\end{rem}

Let $\Delta\in\Hom(\bigwedge^N\cA,\cA)$ be an
$N$-linear skew-symmetric bracket:
$\Delta(a_{\varSigma(1)}$, $\ldots$, $a_{\varSigma(N)})=
(-1)^{\varSigma}\Delta(a_{1},\ldots,a_{N})$ for any rearrangement
$\varSigma\in S_N$.

\begin{define}
The algebra $\cA$ is the \emph{homotopy $N$-Lie} algebra,
or the \emph{Schlessinger\/--\/Stasheff $N$-algebra},
if the $N$-Jacobi identity
\begin{equation}\label{JacobiAsAction}
\Delta[\Delta]=0
\end{equation}
holds.
\end{define}

In coordinates, the $N$-Jacobi identity is
\begin{equation}\label{NaryJacobi}
\sum\limits_{\sigma\in S^N_{2N-1} } (-1)^\sigma\Delta(\Delta(
a_{\sigma(1)},\ldots,
   a_{\sigma(N)}), a_{\sigma({N+1})},\ldots, a_{\sigma(2N-1)}) = 0
\end{equation}
for any $a_j\in\cA$, $1\leq j\leq 2N-1$.
Generally, the number of summands in \eqref{NaryJacobi} is
$\binom{2N-1}{N-1}=\binom{2N-1}{N}$, see \cite{YS2001}.

The Jacobi identity of the type $(N,0,0)$
$\rn{\Delta}{\Delta}=2\Delta[\Delta]=0$
implies Eq.\ \eqref{NaryJacobi} for any even $N$.
If $N$ is odd, then the expression
\begin{equation}\label{RNZero}
\rn{\Delta}{\Delta}=0
\end{equation}
is
trivial and we consider Eq.\ \eqref{JacobiAsAction} separately
from condition \eqref{nkrJacobi}.
In the sequel, we study the Jacobi identity \eqref{NaryJacobi} of
the form \eqref{JacobiAsAction}, where
$\Delta\in\Hom(\bigwedge^N\cA,$ $\cA)$.

\subsection{The Hochschild and the Koszul cohomologies}
The graded Jacobi identity for the
Richardson\/--\/Nijenhuis bracket provides the Hochschild
$\Id_\Delta$-cohomologies on $\Hom(\bigwedge^\ast\cA,$ $\cA)$ for
$\Delta\in\Hom(\bigwedge^k\cA$, $\cA)$, where $k$ is even:
\begin{state}[\cite{Uvas}]\label{HochschildDifferential}
The Richardson\/--\/Nijenhuis bracket
satisfies the graded Jacobi identity
\begin{equation}\label{Jacobi4Schouten}
\rn{\Delta}{\rn{\nabla}{\square}}=
\rn{\rn{\Delta}{\nabla}}{\square} +
   (-1)^{(\Delta-1)(\nabla-1)}
\rn{\nabla}{\rn{\Delta}{\square}}.
\end{equation}
\end{state}
\begin{cor}
Let $k$ be an even natural number
and an operator $\Delta\in\Hom(\bigwedge^k\cA$, $\cA)$
be such that $\rn{\Delta}{\Delta}=0$; then the operator
$
\Id_\Delta\equiv\rn{\Delta}{{\cdot}}
$
is a differential:
$\Id_\Delta^2=0$.
\end{cor}
The cohomologies w.r.t.\ the differential $\Id_\Delta$ are called the
Hochschild $\Id_\Delta$-cohomologies of the space
$\Hom_\Bbbk(\bigwedge^*\cA$, $\cA)$.

\begin{rem}
The approach under study
is closely related with the algebraic Hamiltonian
formalism in the geometry of partial differential equations
(\cite{ClassSymEng, DIPS6-2002}):
a bi-vector $A$ endowes the space of the Hamiltonians with the Lie
algebra structure iff its Schouten bracket with itself satisfies the
equation $\sch{A}{A}=0$; from the Jacobi identity similar to Eq.\
\eqref{Jacobi4Schouten} it follows that the operator
$\Id_A=\sch{A}{{\cdot}}$ is a differential and therefore defines the
Hamiltonian complex whose cohomologies are called the Poisson
cohomologies.
We note that the operator $W^{0,1}=\bun\wedge d/dx$, which is studied
in Section \ref{subsecWronsk}, is the first Hamiltonian structure
for the Korteweg\/--\/de Vries equation.
\end{rem}

Let the bracket $\Delta\in\Hom(\bigwedge^k\cA,\cA)$ satisfy the
homotopy $k$-Jacobi identity $\Delta[\Delta]=0$. By $\dd_\Delta$
denote the linear map
$\dd_\Delta\in\Hom(\bigwedge^r\cA,\bigwedge^{r-k+1}\cA)$ such that
\begin{enumerate}
\item ${\dd_\Delta\Bigr|}_{\bigwedge^r\cA}=0$ if $r<k$;
\item
$\dd_\Delta(a_1\wedge\ldots\wedge a_r) = \sum\limits_{\sigma\in
S_r^k} (-1)^\sigma
\Delta[a_{\sigma(1)},\ldots,a_{\sigma(k)}]\wedge
a_{\sigma(k+1)}\wedge\ldots\wedge a_{\sigma(r)}$ otherwise.
\end{enumerate}
We obtain the Koszul $\dd_\Delta$-cohomologies for the 
algebra $\bigwedge^*\cA=\bigoplus_{r=2}^\infty\bigwedge^r\cA$
over the algebra $\cA$ owing to
\begin{state}[\cite{HanlonWachs}]\label{KoszulDifferential}
The operator $\dd_\Delta\colon\bigwedge^*\cA\to\bigwedge^*\cA$ is a
differential: $\dd_\Delta^2=0$.
\end{state}

By $H^*_\Delta(\cA)$ we denote the Koszul
$\dd_\Delta$-cohomologies w.r.t.\ the differential $\dd_\Delta$.
For $N=2$, the Koszul cohomologies of the Lie algebra of vector fields
on the circumpherence $\mathbb{S}^1$
were obtained in \cite{GelfandFuks68}. For $N\geq2$,
the Koszul $\dd_\Delta$-cohomologies of free algebras were found in
\cite{HanlonWachs}.

\subsection{Examples of the homotopy Lie algebras}
One should notice that algebraic structures \eqref{JacobiAsAction}
have a remarkable geometric motivation to exist.
Namely, we have

\begin{example}[\cite{HanlonWachs}]\label{CrossProduct}\label{HWExample}
Let $\cA=\Bbbk^{r}$ be a $\Bbbk$-linear space of arbitrary
dimension~$r$ and
$\Delta\colon\bigwedge^N\cA\to\cA$ be a
skew-symmetric linear mapping of $\cA$.
If $\dim\cA<2N-1$, then the identity \eqref{NaryJacobi} holds
for~$\Delta$.
\end{example}

\begin{proof}
We maximize the number of summands in \eqref{NaryJacobi}
in order to note its skew-symmetry w.r.t.~the transpositions
$a_j\mapsto a_{\varSigma(j)}$, $\varSigma\in S_{2N-1}$.
%
The l.h.s.~of Jacobi identity (\ref{NaryJacobi}) equals
\begin{equation}\label{AllPermut}
\frac{1}{N!(N-1)!}
\sum\limits_{\sigma\in S_{2N-1} } (-1)^\sigma\Delta(\Delta(
a_{\sigma(1)},\ldots,
   a_{\sigma(N)}), a_{\sigma({N+1})},\ldots, a_{\sigma(2N-1)}),
\end{equation}
where \emph{all} elements $\sigma\in S_{2N-1}$ are taken into
consideration; see \cite{HanlonWachs, MichorVin}.
Expression (\ref{AllPermut}) is skew-symmetric w.r.t.~any
rearrangement $\varSigma$ of the elements $a_j\in\vec{a}$:
\begin{multline}
(-1)^\varSigma
\sum\limits_{\sigma\in S_{2N-1} }\!\!
(-1)^\sigma\Delta(\Delta(a_{(\sigma\circ\varSigma)(1)},\ldots,
a_{(\sigma\circ\varSigma)(N)}),
a_{(\sigma\circ\varSigma)(N+1)},\ldots,
a_{(\sigma\circ\varSigma)(2N-1)}) ={}\\
{}=
\sum\limits_{\sigma\in S_{2N-1} } (-1)^\sigma
   \Delta(\Delta(a_{\sigma(1)},\ldots,
   a_{\sigma(N)}), a_{\sigma({N+1})},\ldots, a_{\sigma(2N-1)}).
\label{SignThrough}
\end{multline}
Consequently, the l.h.s~in (\ref{AllPermut}) is skew-symmetric also
and we obtain a $(2N-1)$-linear skew-symmetric operator acting on the
vector space of smaller dimension.
Hence, if
$\dim\cA<2N-1={}$the number of arguments $\sharp\vec{a}$, then
\eqref{NaryJacobi} holds.
\end{proof}

\noindent$\bullet$
The case $\dim\cA_1=N+1$, $[a_0,\ldots,\widehat{a_j},\ldots,a_N]=
(-1)^{j}\cdot a_j$ is well-known to be the cross-product in
$\Bbbk^{N+1}$. For $\Bbbk=\BBR$ and $N=2$, we have the Lie algebra
$\cA\simeq\mathfrak{so}(3)$.

\noindent$\bullet$
We claim that the algebra $\cA_2$ of dimension $\dim\cA_2=N+1$,
defined by the relations
\begin{equation}\label{Homot4Wronsk}
[a_0,\ldots,\widehat{a_j},\ldots,a_N]=a_{N-j},\quad 0\leq j\leq N,
\end{equation}
admits a representation in the space of polynomials $\Bbbk_N[x]$ such
that its structure $[\ldots]$ is defined
by the Wronskian determinants of scalar fields (smooth
functions of one argument). The algebra $\cA_2$ is considered in
the next subsection.

\subsection{The polynomials}\label{FiniteDimSec}
In this section, we construct finite-dimensional
homotopy $N$-Lie generalizations of the Lie
algebra $\mathfrak{sl}_2(\Bbbk)$. Our starting point is the following
\begin{example}[\textup{\cite{YS2001}}]\label{Starting:sl2}
The polynomials $\Bbbk_2[x]=\{\alpha x^2+\beta x+\gamma\mid
 \alpha,\beta,\gamma\in\Bbbk\}$
of degree $2$ 
form a Lie algebra isomorphic to $\mathfrak{sl}_2(\Bbbk)$.
The latter is generated by three elements $\langle e$, $h$, $f\rangle$
that satisfy the relations
\begin{align}
[h,e]&=2e,  & [h,f]&=-2f,       & [e,f]&=h.\notag \\
\intertext{Consider the basis $\langle1$, $-2x$, $-x^2\rangle$ and choose
the Wronskian determinant for the bracket on $\Bbbk_2[x]$:}
[-2x,1]&=2, & [-2x,-x^2]&=2x^2, & [1,-x^2]&=-2x,\notag \\
\intertext{whence the representation $\rho\colon\mathfrak{sl}_2(\Bbbk)\to
\Bbbk_2[x]$ is}
\rho(e)&=1, & \rho(h)&=-2x, & \text{and}\quad
\rho(f)&=-x^2.\label{RepresSL2}
\end{align}
%
\end{example}

Consider the space $\Bbbk_N[x]\ni a_j$ of polynomials $a_j$ of degree not
greater than $N$; on this space, there is the $N$-linear skew-symmetric
bracket
\begin{equation}\label{NaryBracket}
[a_1, \ldots, a_N]=\W{a_1}{\ldots}{a_N},
\end{equation}
where $W$ denotes the Wronskian determinant.
Since $N$-ary
bracket (\ref{NaryBracket}) is $N$-linear, we consider monomials
$\const\cdot x^k$ only. We choose $\{a^0_j\}=\{x^k\}$ or
$\{a^0_j\}=\{x^k/k!\}$, where $0\leq k\leq N$ and $1\leq j\leq 2N-1$,
for standard basis in $\Bbbk_N[x]$.
The powers
$x^0$, $\ldots$, $x^N$ and $n$ independent variables $x^1$, $\ldots$,
$x^n$ introduced in Section~\ref{nDimWronsk} are never mixed,
and the notation is absolutely clear from the context.
Exact choice of the basis depends on the situation:
the monomials $x^k$ are used to demonstrate the
presence or absence of certain
degrees in $N$-linear bracket \eqref{NaryBracket} and the monomials
$x^k/k!$ are convenient in calculations since they are closed
w.r.t.~the derivations (and the Wronskian determinants as well).
Indeed, we have


\begin{theor}[\textup{\cite{YS2002}}]\label{ClosedDegN}
\label{MonomialsClosed} Let $0\leq k\leq N$\textup{;} then
\label{IsTheMonomial}\label{MonomialCoefficient}
 the relation
\begin{equation}\label{WronskWithConst}
\W{1}{\ldots,\widehat{\frac{x^k}{k!}},\ldots}{\frac{x^N}{N!}}=
\frac{x^{N-k}}{(N-k)!}
\end{equation}
holds.
\end{theor}
\begin{proof}
We have
\begin{equation}\label{WronskDecompose}
\W{1}{\ldots,\widehat{\frac{x^k}{k!}},\ldots}{\frac{x^N}{N!}} =
  \W{1}{\ldots}{\frac{x^{k-1}}{(k-1)!}} \cdot
  \W{x}{\ldots}{\frac{x^{N-k}}{(N-k)!}},
\end{equation}
where the first factor in the r.h.s.~of (\ref{WronskDecompose})
equals $1$ and has the degree $0$. Denote the
second factor, the determinant of the $(N-k)\times(N-k)$ matrix, by
$W_m$, $m\equiv N-k$. We claim that $W_m$ is a monomial:
$\deg W_m=m$, and prove
this fact by induction on $m\equiv N-k$. For $m=1$,
$\deg\det(x)=1=m$. Let
$m>1$; the decomposition of $W_m$ w.r.t.~the last row gives
\begin{equation}\label{DecompositionWRTLastRow}
W_m=\W{x}{{\ldots}}{\frac{x^m}{m!}}=
    x\cdot\W{x}{{\ldots}}{\frac{x^{m-1}}{(m-1)!}} -
    \W{x}{{\ldots},\frac{x^{m-2}}{(m-2)!}}{\frac{x^m}{m!}},
\end{equation}
where the degree of the first Wronskian in r.h.s.~of
(\ref{DecompositionWRTLastRow}) is $m-1$ by the inductive assumption.
Again, decompose the second Wronskian in r.h.s.~of
(\ref{DecompositionWRTLastRow}) w.r.t.~the last row and proceed
iteratively by using the induction hypothesis. We obtain
the recurrence relation
\begin{equation}\label{RecurrentWronsk}
W_m=\sum\limits_{l=1}^{m-1}
W_{m-l}\cdot (-1)^{l+1} \frac{x^l}{l!} -
(-1)^m\,\frac{x^m}{m!},\qquad m\geq1,
\end{equation}
whence $\deg W_m=m$.
We see that the initial Wronskian (\ref{WronskDecompose}) is a
monomial itself of degree $m=N-k$ with yet unknown coefficient.

Now we calculate the coefficient $W_m(x)/x^m\in\Bbbk$ in the Wronskian
determinant \eqref{WronskDecompose}. Consider the generating function
\begin{equation}\label{GenFunction}
f(x)\equiv\sum\limits_{m=1}^\infty W_m(x)
\end{equation}
such that
$$
W_m(x)=\frac{x^m}{m!}\,\frac{d^mf}{dx^m}(0),\qquad 1\leq m\in\BBN.
$$
Recall that $\exp(x)\equiv\sum_{m=0}^\infty x^m/m!$;
treating \eqref{GenFunction} as the formal sum of equations
\eqref{RecurrentWronsk}, we have
$$
f(x)=f(x)\cdot(\exp(-x)-1)-\exp(-x)+1,
$$
whence
\begin{equation}\label{ComputeGenFunction}
f(x)=\exp(x)-1.
\end{equation}
Hence the required coefficient equals $1/m!$.
The proof is complete.
\end{proof}

We have shown that the polynomials $\Bbbk_N[x]$ are closed w.r.t.\ the
Wronskian determinant, and we know that any $N$-linear skew-symmetric
bracket $\Delta$ on $\Bbbk^{N+1}$ satisfies $\Delta[\Delta]=0$.
Therefore,
\label{TheorDegN}
the statement that the
polynomials $\Bbbk_N[x]$ of degree not greater than $N$ form the
homotopy $N$-Lie algebra with $N$-linear skew-symmetric bracket
\eqref{NaryBracket} for any integer $N\geq2$ is quite obvious.
Nevertheless, in the sequel we show that the Wronskian
$W^{0,1,\ldots,N-1}\in\Hom(\bigwedge^N\Bbbk_N[x],\Bbbk_N[x])$ is the
restriction of a \emph{nontrivial}
homotopy $N$-Lie bracket that lies in
$\Hom(\bigwedge^N\Bbbk[[x]],\Bbbk[[x]])$.
Also, the dimension $n$ of the base $\Bbbk\equiv\Bbbk^1\ni x$ equals
$1$.  In Section \ref{nDimWronsk},
we generalize the concept to the case $x\in\Bbbk^n$,
where integer $n\geq1$ is arbitrary.

\section{The associative homotopy Lie algebras}\label{SecAssocAlg}
\noindent%
Another natural example of the homotopy Lie algebras is given by
\begin{state}[\cite{DzhFAP, HanlonWachs}]\label{Associative}
Let $\cA$ be an associative algebra and let $N$ be
even\textup{;} by definition, put\footnote{Note that the
permutations $\sigma\in S_N$ provide the \emph{direct} left action on
$\bigotimes^N\!\cA\!$, contrary to the inverse action in
\cite[\S{}II.2.6]{Kassel}.
By definition, $\sigma(j)$ is the index of the object in an
initial ordered set placed onto $j$th position after the left action
of a permutation $\sigma$.}
\begin{equation}\label{AssocBracketHW}
[a_1,\ldots,a_N]\stackrel{\mathrm{def}}{=}
\sum\limits_{\sigma\in S_{N}} (-1)^\sigma\cdot
a_{\sigma(1)}\circ\dots\circ a_{\sigma(N)}.
\end{equation}
Then $\cA$ is a homotopy Lie algebra w.r.t.\ this bracket.
\end{state}
\begin{proof}[Proof \textup{(\cite{HanlonWachs})}.]
The crucial idea is using \eqref{AllPermut} and \eqref{SignThrough}. Let
$a_1$, $\ldots$, $a_{2N-1}$ lie in $\cA$ and $\sigma\in S_{2N-1}$ be a
permutation. In order to compute the coefficient of
$a_{\sigma(1)}\circ\dots\circ a_{\sigma(2N-1)}$
in \eqref{NaryJacobi} and prove
it to be trivial, it is enough to do that for
$\alpha=a_1\circ\dots\circ a_{2N-1}$ in \eqref{NaryJacobi} due to
\eqref{SignThrough} and \eqref{AllPermut}, successively.

Now we use the assumption $N\equiv0\mod2$. The product $\alpha$ is met
$N$ times in \eqref{NaryJacobi} in the summands $\beta_j$, $1\leq j\leq
N$:
\begin{equation}\label{HWCycle}
\beta_j=(-1)^{N(j-1)}\,
[[a_j,\ldots,a_{N+j-1}],a_{1},\ldots,a_{j-1},a_{N+j},\ldots,
   a_{2N-1}].
\end{equation}
The coefficient of $\alpha$ in $\beta_j$ equals $(-1)^{j-1}$ and hence
the coefficient of $\alpha$ in \eqref{NaryJacobi} is
$$
\sum_{j=1}^N(-1)^{j-1}=0.
$$
The proof is complete.
\end{proof}

From the proof of Proposition \ref{Associative} we see that the main
obstacle for bracket \eqref{AssocBracketHW} to provide the homotopy
Lie algebra structures for odd $N$s are the signs within
\eqref{Action}, \eqref{AssocBracketHW},  and in the
Richardson\/--\/Nijenhuis bracket \eqref{RNBracket} that defines the
Jacobi identity as the cohomological conditions $\Id_\Delta^2=0$,
see Proposition \ref{HochschildDifferential}. Namely, we have
\begin{state}[\cite{DzhFAP}]\label{AssocObstaclesOdd}
Let the subscript $i$ at the bracket's \eqref{AssocBracketHW}
symbol $\Delta_i$ denote its number of arguments\textup{:}
$\Delta_i\in\Hom_\Bbbk(\bigwedge^i\cA,\cA)$, and let $k$ and $\ell$ be
arbitrary integers. Then the following identities hold\textup{:}
\begin{subequations}
\begin{align}
\Delta_{2k}[\Delta_{2\ell}]&=0,\label{BothEven}\\
\Delta_{2k+1}[\Delta_{2\ell}]&=\Delta_{2k+2\ell},\label{InnerEven}\\
\Delta_k[\Delta_{2\ell+1}]&=k\cdot\Delta_{2\ell+k}.\label{InnerOdd}
\end{align}
\end{subequations}
\end{state}
\begin{proof}
The proof of \eqref{BothEven} repeats the reasoning in \eqref{HWCycle}
literally. For \eqref{InnerEven}, we note that the last
summand $\beta_{2k+1}$ is not compensated. For \eqref{InnerOdd},
the summand $\alpha=a_1\circ\dots\circ a_{2\ell+k}$ acquires the
coefficient
$$
\sum_{j=1}^k(-1)^{(2\ell+1)(j-1)}\cdot(-1)^{j-1}=k.
$$
This completes the proof.
\end{proof}

\subsection{On representations in differential operators}
In this section, the field $\Bbbk$ is the complex field $\BBC$:
$\Bbbk=\BBC$, and $z$ is the holomorphic coordinate in $\BBC$.

By $\cO(\BBC)$ we denote the algebra of the Laurent series over $\BBC$.
Consider the associative algebra $\Diff_*(\cO(\BBC),\cO(\BBC))$ of
holomorphic differential operators
\begin{equation}\label{Htransform}
\nabla_{\vec w}=\sum\limits_{j=0}^p w_j(z)\cdot\dd^j.
\end{equation}
We claim that for any $p$ the algebra $\Diff_*(\cO(\BBC),\cO(\BBC))$
is a homotopy Lie algebra that
contains a homotopy $2p$\/-\/Lie subalgebra
defined by a skew\/--\/symmetric bracket of $2p$ arguments.
Recall that for $N=2$
vector fields compose a subalgebra in the space of
differential operators of arbitrary orders.

Let $a_j\in\Diff_*(\BBC)$ be $a_j=w_j(z)\,\dd^{k_j}$
for $1\leq j\leq N$.
Similar to~\eqref{AssocBracketHW}, put
\begin{equation}\label{AssocBracket}
[w_1\cdot\dd^{k_1},\ldots,w_N\cdot\dd^{k_N}]\stackrel{\mathrm{def}}{=}
\sum\limits_{\sigma\in S_N}
(-1)^\sigma\,w_{\sigma(1)}\,\dd^{k_{\sigma(1)}}\circ\dots\circ
w_{\sigma(N)}\cdot\dd^{k_{\sigma(N)}}.
\end{equation}
This bracket is $N$-linear over $\BBC$ and is skew-symmetric w.r.t.\
permutations of its arguments.

First, we count derivatives: Consider the special case $k_j\equiv
p=\const$ for all $j$ and solve the equation
\begin{equation}\label{Balance}
Np=\frac{N(N-1)}{2}+p
\end{equation}
for $p$: $p={N}/{2}$; note that $N(N-1)/2=|W^{0,1,\ldots,N-1}|$.

Further on, we restrict ourselves to the case $N\equiv0\mod2$; it turns
out that for odd $N$s we need to consider half\/-\/integer powers of
the derivation $\dd$: $\dd^0$, $\dd^{1/2}$, $\dd$, $\dd^{3/2}$,
$\ldots$. 

\begin{theor}\label{OnlyWronsk}
Let $N$ be even and $w_j\in\cO(\BBC)$ for $0\leq j\leq N$\textup{;}
put $p=N/2$. Then we have
$$
[w_1\,\dd^p,\ldots,w_N\,\dd^p]=W^{0,1,\ldots,N-1}(w_1,\ldots,w_N)\cdot
\dd^p,
$$
where $W^{0,1,\ldots,N-1}=\bun\wedge\dd\wedge\ldots\wedge\dd^{N-1}$ is
the Wronskian determinant.
\end{theor}

\begin{proof}
Permutations of arguments in the r.h.s.\ of \eqref{AssocBracket} are
reduced to permutations of $w_j$s since $k_j\equiv p$. Let $\sigma\in
S_N$ be a permutation and $\vec\jmath\in\BBZ^N\cap[0,Np]^N$ be a vector
in the integral lattice. Suppose that the r.h.s.\ in
\eqref{AssocBracket} is expanded from left to right and all possible
derivation combinations
$$
S_{\sigma,\vec\jmath}\stackrel{\mathrm{def}}{=}
(-1)^\sigma\,\dd^{j_1}(w_{\sigma(1)})\cdot\ldots\cdot
\dd^{j_N}(w_{\sigma(N)})
$$
are obtained; we note that \emph{not} all vectors
$\vec\jmath\in\BBZ^N\cap[0,Np]^N$ can be realized: at least,
$|\vec\jmath|\leq Np$. Still, the set
$J=\{\vec\jmath\}\subset\BBZ^N\cap[0,Np]^N$ does not depend on
$\sigma$. Assume there is a summand such that two functions $w_a$ and
$w_b$ acquire equal numbers of derivations for some combination
$\vec\jmath\in J$. Then, for the same combination $\vec\jmath$
and the transposition $\tau_{ab}$, there is
the permutation $\tau_{ab}\circ\sigma$ such that the order of
$w_a$ and $w_b$ is reversed and $S_{\sigma,\vec\jmath}+
S_{\tau_{ab}\circ\sigma,\vec\jmath}=0$. Consequently, only the Wronskian
remains at $\dd^p$ owing to Eq.~\eqref{Balance}.
\end{proof}

Theorem \ref{OnlyWronsk} is a generalization of a perfectly familiar
fact: the commutator of two vector fields is a vector field.
We emphasize that Theorem \ref{OnlyWronsk} forbids the naive approach
that combines $N$ vector fields (e.g., symmetries of a PDE) in an
attempt to obtain some vector field again. 

\begin{rem}\label{Unfortune}
Unfortunately, for arbitrary operators \eqref{Htransform} of
order $p=N/2$ this mechanism of compensations does not work. Indeed,
suppose that the powers $k_j\leq p$ are arbitrary; then the sets
$J\subset\BBZ^N\cap[0,Np]^N$ do \emph{depend on} $\sigma$, and
generally $J(\sigma)\neq J(\tau_{ab}\circ\sigma)$
if two functions $w_a$ and $w_b$ are differentiated w.r.t.\ $z$ equal
number of times in a summand $S_{\sigma,\vec\jmath(\sigma)}$. Of
course, we can obtain the Wronskian determinant at some suitable power
of $\dd$, but there can be much more summands, even at $\dd^\ell$ for
$\ell\geq p$.
The same difficulty occurs for the lower bound
$k_j\geq p$, when we consider formal
differential operators $\nabla=\sum_{j\geqslant p} w_j(z)\cdot\dd^j$.

Nevertheless, for an arbitrary integer $p'\geq(N-1)/2$ we have
$$
[w_1\,\dd^{p'},\ldots,w_N\,\dd^{p'}]=
W^{0,1,\ldots,N-1}(w_1,\ldots,w_N)\cdot
\dd^{Np'-N(N-1)/2}.
$$
\end{rem}

For various pairs $(N,p)\in\BBN\times\BBN$, one can deduce many
extravagant phenomena. In \cite{DzhCommutators}, the following
proposition is proved: for $p=1$, vector fields $\ID(M^n)$ on a smooth
$n$\/-\/dimensional manifold $M^n$ are closed w.r.t\ $N$\/-\/ary
bracket \eqref{AssocBracket} and form the homotopy $N$\/-\/Lie
algebra if $N=n^2+2n-2$.

As a corollary to Proposition \ref{Associative}, we have

\begin{theor}\label{ClosedWithWronskians}
Let $N$ be even\textup{;}
consider the $\cO(\BBC)$-module
$W_{N/2}\stackrel{\mathrm{def}}{=}{\mathrm{span}}_{\BBC}
\langle w(z)$ $\dd^{N/2}\rangle$ of holomorphic operators of
order~$N/2$. Then $W_{N/2}$ is endowed with the homotopy $N$-Lie
algebra structure w.r.t.\ bracket~\eqref{AssocBracket}.
\end{theor}

Nevertheless, the difficulties in complete description of the
r.h.s.\ in \eqref{AssocBracket} do not influence upon our ability to
observe the homotopy $N$-Lie structure on the associative algebra of
operators \eqref{Htransform}.
As a reformulation to Proposition \ref{Associative}
on page \pageref{Associative} we obtain
\begin{theor}\label{EvenNsTh}
Let $N$ be even, then differential operators
\eqref{Htransform} of
arbitrary orders compose the homotopy $N$-Lie algebra w.r.t.\
bracket \eqref{AssocBracket}.
\end{theor}

Indeed, the differential operators generate the associative algebra
$\Diff_*(\cO(\BBC),\cO(\BBC))$.

\subsection{On the Wronskian determinants}\label{subsecWronsk}
We start with
\begin{state}[\cite{DzhFr}]\label{J2Zero}
Let $k$ and $l$ be positive integers, then the identity
$$
W^{0,1,\ldots,k}[W^{0,1,\ldots,l}]=0
$$
holds.
\end{state}
\begin{rem}\label{ForwardProofRem}
Actually, a slight modification of Theorem \ref{ClosedWithWronskians}
combined with Proposition \ref{AssocObstaclesOdd}
give a nice and compact proof of Proposition \ref{J2Zero} in the case
when the numbers $k$ and $l$ of the arguments are even.
In the next section, we generalize Proposition \ref{J2Zero} to the
Wronskians $D^{\vec\sigma}$ w.r.t.\ several independent variables
$x^1$, $\ldots$, $x^n$, and, in particular, obtain its proof
for arbitrary naturals $k$ and $l$ in the case $n=1$.
\end{rem}

From Proposition \ref{J2Zero} we also obtain
\begin{theor}\label{RNBracketZero}
Let $k$ and $l$ be positive integers, then the relation
$$
\rn{W^{0,1,\ldots,k}}{W^{0,1,\ldots,l}}=0
$$
holds.
\end{theor}

\begin{cor}
The $\Id_W$-cohomologies of the space of Wronskians are isomorphic to
this space itself:
$H_{\Id_W}^*=\mathrm{span}_\Bbbk\langle W^{0,\ldots,l}$,
$l\geq1\rangle$, since the differential $\Id_{W^{0,\ldots,k}}=
\rn{W^{0,\ldots,k}}{{\cdot}}$ is trivial.
\end{cor}

\label{OneDimWronskSmooth}

Now we study the homotopy
generalizations of the Witt algebra
(the Virasoro algebra with zero central charge)
which is defined by the relations
$[a_i, a_j]=(j-i)\,a_{i+j}$.
Taking into account all our observations
on the Wronskians, we set $N=2$ and consider the polynomial generators
$a_i=x^{i+1}$,
where $x\in\Bbbk$ and $i\in\BBZ$. For $N\geq2$ and
the Wronskian determinant
$W^{0,1,\ldots,N-1}$, we consider the relations
\begin{equation}\label{WittGeneral}
[a_{i_1},\ldots,a_{i_N}]=\Omega(i_1,\ldots,i_N)\,a_{i_1+\cdots+i_N},
\end{equation}
where the structural constants $\Omega(i_1,\ldots,i_N)$ are
skew-symmetric w.r.t.\ their arguments. Here we use the representation
$a_i=x^{i+N/2}$. We claim that the function $\Omega$ is the Vandermonde
determinant.

\begin{theor}\label{VanderAsCoeff}
Let $\nu_1$, ${\ldots}$, $\nu_N\in\Bbbk$ be constants and set
$\nu=\sum_{i=1}^N\nu_i$\textup{;} then the equality
\begin{equation}\label{Vander}
W^{0,1,\ldots,N-1}(x^{\nu_1},\ldots,x^{\nu_N})=
\prod\limits_{1\leq i<j\leq N}(\nu_j-\nu_i)\cdot
x^{\nu-{N(N-1)}/{2} }
\end{equation}
holds, i.e., the Wronskian determinant of the
monomials is a monomial
itself and its coefficient is the Vandermonde determinant.
\end{theor}
\begin{proof}
Consider the determinant \eqref{Vander}:
$A=\det\|a_{ij}\,x^{\nu_j-i+1}\|$. From $j$th column take the monomial
$x^{\nu_j-N+1}$ outside the determinant:
$$
A=x^{\nu-N(N-1)}\cdot\det\|a_{ij}\,x^{N-i}\|;
$$
all rows acquire common degrees in $x$: $\deg$(any element in $i$th
row) ${}=N-i$. From $i$th row take this common factor $x^{N-i}$ outside
the determinant:
$$
A=x^{\nu-N(N-1)/2}\cdot\det\|a_{ij}\|,
$$
where the coefficients $a_{ij}$ originate from the initial derivations
in a very special way:
for any $i$ such that $2\leq i\leq N$, we have
$$
a_{1j}=1\qquad\text{and}\qquad
a_{ij}=(\nu_j\underline{{}-i+2})\cdot a_{i-1,j}
\quad\text{for $1<i\leq N$.}
$$
The underlined summand does not depend on $j$, and hence for any $k=N$,
${\ldots}$, $2$ the determinant $\det{}\|a_{ij}\|$ can be splitted in
the sum:
\begin{multline*}
\det\|a_{ij}\|=\det\|a'_{kj}=\nu_j\cdot a_{k-1,j};\quad
   a'_{ij}=a_{ij} \text{ if }i\neq k\| +
 {}\\
    {}+  \det\|a''_{kj}=(2-i)\cdot a_{k-1,j};\quad
   a''_{ij}=a_{ij} \text{ if }i\neq k\|,
\end{multline*}
where the last determinant is trivial.

Solving the recurrence relation, 
we obtain
$$
\det\|a_{ij}\|=\det\|\nu_j^{i-1}\|=\prod\limits_{1\leq k<l\leq N}
(\nu_l-\nu_k).
$$
This completes the proof.
\end{proof}

\begin{rem}
We have calculated the structural constants in \eqref{WittGeneral}
by using
another basis $a_i'=x^i$ such that the resulting degree is not
$\sum_{k=1}^N\deg a_k'$. Nevertheless, the result is correct since we
use the translation invariance of the Vandermonde determinant:
$$
\Omega(i_1,\ldots,i_N)=\Omega(i_1+\frac{N}{2},\ldots,i_N+\frac{N}{2}),
$$
and therefore all reasonings are preserved.
\end{rem}


Now we recall the behaviour of bracket \eqref{NaryBracket} w.r.t.\ a
change of coordinates $y=y(x)$.
\begin{theor}\label{ConformalWeightTh}
Let $\phi^i(y)$ be smooth functions for $1\leqslant i\leqslant N$, i.e.,
$\phi^i$ is a scalar field of the conformal weight $0$, such that $\phi^i$
is transformed by the rule $\phi^i(y)\mapsto\phi^i(y(x))$ under a change
$y=y(x)$.
Then the relation
$$
{\det\left\|\frac{d^j\phi^i}{dx^j}\right\|}_{
 \begin{array}{rcl}
    i&\!=\!&1,\ldots,N\\
    j&\!=\!&0,\ldots,N-1\\
    \phi^i&\!=\!&\phi^i(y(x))
 \end{array}
 }
= {\left(\frac{dy}{dx}\right)}^{\Delta(N)}\,
{\det\left\|\frac{d^j\phi^i}{dy^j}\right\|}_{
 \begin{array}{rcl}
    \phi^i&\!=\!&\phi^i(y)\\
    y&\!=\!&y(x)
 \end{array}
 }
$$
holds, where the conformal weight $\Delta(N)$ for the Wronskian
determinant of $N$ scalar fields $\phi^i$ of weight $0$ is
$\Delta(N)={N(N-1)}/{2}$.
\end{theor}
\begin{proof}
Consider a function $\phi^i(y(x))$ and apply the total derivative
$D_x^j$ by using the chain rule. The result is
$$
\frac{d^j\phi^i}{dy^j}\cdot{\left(\frac{dy}{dx}\right)}^j +
\text{terms of lower order derivatives }\frac{d^{j'}}{dy^{j'}},\quad
j'<j.
$$
These lower order terms differ from the leading terms in
$D_x^{j'}\,\phi^i(y(x))$, $0\leq j'<j$, by the factors common for all
$i$ and thus they produce no effect since a determinant with
coinciding (or proportional) lines equals zero.
From $i$th row of the Wronskian we extract $(i-1)$th power of $dy/dx$,
their total number being $N(N-1)/2$. This is the conformal weight
by definition.
\end{proof}

We see that the Wronskian determinant of $N$ functions is \emph{not} a
function itself: the objects we are dealing with are the higher
order differential operators, and the functions are their coefficients
w.r.t.\ the basis $\langle\bun$, $\dd$, $\ldots\rangle$.

Theorem \ref{ConformalWeightTh} can be extended to the case $n\geq1$:
\textbf{\textit{x}}${}=(x^1$, $\ldots$, $x^n)$. We also see
that the statement is generally not true if the generalized Wronskian
is $\dd^{\sigma_1}\wedge\ldots\wedge\dd^{\sigma_N}\neq\const\cdot
\mathbf{1}\wedge\dd\wedge\ldots\wedge\dd^{N-1}$.

\subsection{Applications in the $W$\/--\/geometry}\label{WGeometry}
We note that the concept of the homotopy Lie structures for
differential operators \eqref{Htransform} has a nice application in
the $W$-geometry.
The $A_\ell$-$W$-geometry \cite{GervaisLetter}--\cite{GervaisENS}
is the geometry of complex curves
$\varSigma$: $\dim_\BBC(\varSigma)=1$,
$\dim_{\overline{\BBC}}(\varSigma)=1$,
chirally embedded into the K\"ahler manifold
$\BBC\mathbb{P}^\ell$; further on, $f^A(z)$ and $\bar f^{\bar
A}(\bar z)$ are the embedding functions, $0\leq A\leq\ell$ and
$0\leq\bar A\leq\ell$.
The compatibility conditions of these embeddings are the
Toda equations (\cite{LeznovSavelievBV}),
associated with the semisimple $A_\ell$\/-\/type Lie algebras.

\begin{define}[{\cite[\S3.2]{GervaisClassic}}]
A general infinitesimal $W$\/--\/\emph{transformation} $\delta_W$
of such a curve is a
change of the embedding functions $f^A$, $\bar f^{\bar A}$ of the form
\begin{equation}\label{Wtransform}
\delta_W\,f^A(z)=\sum\limits_{j=0}^k w_j(z)\,\dd^jf^A(z),\qquad
\delta_W\,\bar f^{\bar A}(\bar z)=\sum\limits_{j=0}^k
\bar w_j(\bar z)\,\bar\dd^j\bar f^{\bar A}(\bar z),
\end{equation}
where $w_j\in\cO(\BBC)$ and $\bar w_j$are holomorphic and
anti\/-\/holomorphic functions, respectively.
\end{define}

We see that a $W$-transformation is uniquely defined by the higher
order differential operator $\sum_jw_j(z)\cdot\dd^j\in
\Diff_*(\cO(\BBC), \cO(\BBC))$.
So, 
by using Theorem \ref{EvenNsTh} we conclude that
%
the higher\/-\/order $W$\/-\/transformations
compose the homotopy $N$-Lie algebras for even natural~$N$s.

\section{The Wronskians in multidimensions: $n\geq1$}\label{nDimWronsk}%
\label{SecnDimBase}
In this section, we generalize the concept of the
Wronskian determinants to the multidimensional case of the base
$\bx\in\BBR^n$; here we assume $\Bbbk=\BBR$.
Further on, we consider the $k$th order jets $J^k(n,1)$ over
the bundle $\pi\colon\BBR^n\times\BBR\to\BBR^n$, where the base
is the Euclidean space of
dimension $n\geq1$ and the algebra $\cA$ is the associative
commutative algebra $C^\infty(\BBR^n)$ of smooth functions.

In order to construct a natural $n$-dimensional base generalization of
the Wronskians, we pass to the geometrical standpoint
(\cite{ClassSymEng}) and make an experimental observation first.

Let $\cF(\pi)$ be the algebra of smooth
functions $C^\infty(J^\infty(\pi))$ and consider
the $\cF(\pi)$-module $\varkappa(\pi)$ of
evolutionary vector fields\footnote{To denote evolutionary vector
fields, we use the Cyrillic letter $\cEv$, which is pronounced like
``e" in ``ten".}
$\cEv_a=\sum_{j,\sigma}D_\sigma(a^j)\cdot\dd/\dd
u_\sigma^j$, where $a^j\in\cF(\pi)$.
To each Cartan $N$-form
$\omega\in\cC^N\Lambda(\pi)$ we assign the operator
$\nabla_\omega\in\CDiff^{\mathrm{alt}}_{(N)}(\varkappa(\pi)$,
$\cF(\pi))$ by the rule
\begin{equation}\label{FormOperatorIso}
\nabla_\omega(a_1,\ldots,a_N)=\cEv_{a_N}\inner(\ldots
(\cEv_{a_1}\inner\omega)\ldots),
\end{equation}
where $a_i\in\varkappa(\pi)$.

\begin{state}[{\cite[\textup{Chapter 5}]{ClassSymEng}}]
Correspondence \eqref{FormOperatorIso} is the isomorphism of the
$\cF(\pi)$-modules\textup{:}
$$
\cC^N\Lambda(\pi)\simeq
\CDiff^{\mathrm{alt}}_{(N)}(\varkappa(\pi),\cF(\pi)).
$$
\end{state}
Further on, we use the notation $\omega(\cEv_{a_1},\ldots,\cEv_{a_N})
\stackrel{\mathrm{def}}{=} \cEv_{a_N}\inner(\ldots
(\cEv_{a_1}\inner\omega)\ldots)$, where $\omega\in\cC^N\Lambda(\pi)$
and $a_i\in\varkappa(\pi)$.

\begin{rem}[\cite{YS2001}]
Consider the infinite jets $J^\infty(\pi)$ over the bundle
$\pi\colon\BBR\times\BBR\to\BBR$.
Let $x\in\BBR$ be the independent
base variable, $u$ be the dependent fiber variable, $D_x$ be the total
derivative w.r.t.\ $x$, and $u^{(k)}\equiv D_x^k\,u$ be the coordinates
in $J^\infty(\pi)$ for any $k\geq0$. By $\Id_\cC $ we denote the Cartan
differential, $\Id_\cC \colon C^\infty(J^\infty(\pi))\to
\cC\Lambda(J^{\infty}(\pi))$, that maps $u^{(k)}\mapsto
\Id u^{(k)}-D_xu^{(k)}\,dx$.
The Wronskian determinants
\eqref{NaryBracket} can be interpreted as action of the 
$N$-forms
$\Id_\cC u\wedge\ldots\wedge \Id_\cC u^{(N-1)}\in
\cC\Lambda^*(J^{N-1}(\pi))\subset\cC\Lambda^*(J^{\infty}(\pi))$
upon the evolutionary vector fields
$\cEv_{a_j}\equiv\sum_{k=0}^\infty D_x^k(a_j)\,\dd/\dd u^{(k)}$:
\begin{align*}
[a_1,a_2] &= \Id u\wedge\Id(u')\,(\cEv_{a_1},\cEv_{a_2}),\\
[a_1,a_2,a_3] &=\Id u\wedge\Id(u')\wedge
     \Id(u'')\,(\cEv_{a_1},\cEv_{a_2},\cEv_{a_3}),\qquad\text{etc.}
\end{align*}
for any $a_j\in C^\infty(\BBR)$. We emphasize that $a_j\in
C^\infty(\BBR)\subset\varkappa(\pi)$, i.e., we restrict ourselves to
a submodule of $\varkappa(\pi)$ generated by
functions on the base~$M$.
\end{rem}

\begin{rem}
Consider the ternary bracket $\square\in\Hom(\bigwedge^3
C^\infty(\BBR^2), C^\infty(\BBR^2))$:
$$
\square(a_1\wedge a_2\wedge a_3)
=\Id_\cC u\wedge \Id_\cC u_x\wedge \Id_\cC u_y
(\cEv_{a_1},\cEv_{a_2},\cEv_{a_3})=
\begin{vmatrix}
a_1 & a_2 & a_3\\
D_x(a_1) & D_x(a_2) & D_x(a_3)\\
D_y(a_1) & D_y(a_2) & D_y(a_3)
\end{vmatrix}.
$$
For the bracket $\square$, the homotopy ternary Jacobi identity
$\square[\square]=0$ of the form \eqref{JacobiAsAction}
holds. We prove this fact by direct calculations using the \textsf{Jet}
software \cite{Jet97}.
\end{rem}

\begin{state}[\cite{ClassSymEng}]
The dimension of the  vertical part
$J^k(n,1)/\BBR^n$ of the jet space $J^k(n,1)$ is
\begin{equation}\label{DimJets}
\dim\frac{J^k(n,1)}{\BBR^n}=
\dim J^k(n,1)-n=\sum_{i=0}^k\binom{n+i-1}{n-1}=\binom{n+k}{n}.
\end{equation}
\end{state}

We also note that this dimension $N\equiv\binom{n+k}{n}$
is such that the inequality
$$
\dim J^{k_1+k_2}(n,1)-n-1\geq
  \dim J^{k_1}(n,1)+\dim J^{k_2}(n,1)-2(n+1)
$$
is valid for any $k_1$ and $k_2$;
in what follows, we need to substract the dimension $\dim
J^0(n,1)=n+1$ in order to deal with non-trivial multiindexes
$\sigma$ such that $|\sigma|>0$.

Choose arbitrary positive integers $n$ and $k$; then
$N=\binom{n+k}{n}$ is the dimension $\dim(J^k(n,$
$1)/\BBR^n)$.
Let $\cA=C^\infty(\BBR^n)$ be the algebra of smooth functions
$a_j\in\cA$, $1\leq j\leq N$. Now we define the $N$-linear
skew-symmetric bracket $\square_k\in\Hom(\bigwedge^N\cA,\cA)$: we put
\begin{equation}\label{NaryVolumeBracket}
\square_k(a_1,\ldots,a_N)=\bigwedge_{l=0}^k\Bigl(
\bigwedge_{|\sigma|=l}\Id_\cC \cdot D_\sigma u\Bigr)\,
(\cEv_{a_1},\ldots,\cEv_{a_N}).
\end{equation}
In coordinates, this bracket is
$\square_k(a_1,\ldots,a_N)=\det\|D_{\sigma_i}(a_j)\|$,
where $\sigma_i=(\sigma^i_1,\ldots,\sigma^i_n)$ runs through
all multi\/-\/indexes such that $u_{\sigma_i}$ is a coordinate on the
$k$th jet space $J^k(n,1)$ of the fibre bundle
$\BBR\times\BBR^n\to\BBR$.

We claim that the $N$-linear skew-symmetric bracket
$\square_k\in\Hom(\bigwedge^N\cA,\cA)$ defined in
\eqref{NaryVolumeBracket} satisfies the homotopy
$N$-Lie Jacobi identity
\begin{equation}\label{nDimJacobiAsAction}
\square_k[\square_k]=0.
\end{equation}
To prove that, we establish a substantially more general

\begin{theor}\label{nDimBaseJacobi}
Let $N_\nel=\binom{n+k_\nel}{n}$ and $N_\out=\binom{n+k_\out}{n}$ be
the dimensions given by Eq.\ \eqref{DimJets} for some natural $k_\nel$
and $k_\out$\textup{;} by $\square_\nel$ and $\square_\out$ denote the
multilinear skew\/--\/symmetric brackets defined in Eq.\
\eqref{NaryVolumeBracket}. Then the equality
\[
\square_{k_\out}[\square_{k_\nel}]=0
\]
holds.
\end{theor}
\begin{proof}
Without loss of generality we assume that $k_\nel\geq k_\out$,
otherwise one
has to transpose the subscripts '$\nel$' and '$\out$'
in Eq.~\eqref{MinimalNontrivial}.

In contrast with the reasoning in Section \ref{OneDimWronskSmooth},
we deal with
$D^{\vec\sigma}=D_{\sigma_1}\wedge\ldots\wedge D_{\sigma_N}$,
where $\sigma_j$ is a multiindex $(\sharp x^1,\ldots,\sharp
x^n)\in\BBZ^n_{+}$ for any $j$, $1\leq j\leq N=\binom{n+k}{n}$.
Define the norm
$|D^{\vec\sigma}|=|\vec\sigma|=\sum_{j=1}^{N}|\sigma_j|$; we see that
$|\square_{k_\out}[\square_{k_\nel}]|=
|\square_{k_\nel}|+|\square_{k_\out}|$.

Now we note that the non-trivial skew-symmetric
$(N_\nel+N_\out-1)$-linear
bracket $\square_{\min{}}\in\Hom(\bigwedge^{N_\nel+N_\out-1}\cA,\cA)$
with the minimal norm is
\begin{equation}\label{MinimalNontrivial}
\square_{\min{}}=\square_{k_\nel}\wedge\Bigl(
\sum\limits_{\bar\jmath\in\Lambda^{N_\out-1}
\left(J^{k_\nel+k_\out}(n,1)/J^{k_\nel}(n,1)\right)}
\const(\bar\jmath)\cdot D^{\sigma_{\bar\jmath}}
\Bigr),
\end{equation}
where $\const(\bar\jmath)\in\BBR$ are some constant coefficients.

We claim that
$|\square_{\min{}}|>|\square_{k_\nel}+\square_{k_\out}|$,
and thence
$\square_{k_\out}[\square_{k_\nel}]=0$. Indeed, consider the r.h.s.\ in
\eqref{MinimalNontrivial} and note that
$|\Delta\wedge\nabla|=|\Delta|+|\nabla|$, see page
\pageref{exterior}
for definition of the wedge product ${\wedge}$ in this case.
The set of $N_\out$ different
derivatives in $\square_{k_\out}$ admits the canonical splitting:
$$
\square_{k_\out}=D^{\vec\tau^\out}=\mathbf{1}\wedge
\underbrace{D_{\tau_2^\out}\wedge\ldots\wedge
D_{\tau^\out_{N_\out}}}_{\text{$N_\out-1$ factors}},
$$
where $\vec\tau^\out$
contains all multiindexes in $J^{k_\out}(n,1)$, and those
underbraced derivatives are in bijective correspondence with
$N_\out-1$
different derivatives within any summand in the second wedge factor
of \eqref{MinimalNontrivial}
(there is the correspondence owing to the equal numbers of elements).
Still,
$$
1\leq|D_{\tau^\out_i}|=|\tau^\out_i|\leq k_\out <
k_\nel+1\leq|\sigma_{\bar\jmath,i}|=|D^{\sigma_{\bar\jmath,i}}| \leq
k_\nel+k_\out\qquad \forall i\neq1,\quad\forall\bar\jmath.
$$
Indeed, if a multiindex $\sigma_{\bar\jmath,i}$ is such that
$u_{\sigma_{\bar\jmath,i}}$ is a coordinate on the jet space's part
$J^{k_\nel+k_\out}(n,1)/J^{k_\nel}(n,1)$ with the higher order
derivatives only, then
$\sigma_{\bar\jmath,i}$ is \emph{longer} than any multiindex
$\tau^\out_i$ such that $u_{\tau^\out_i}$ is a coordinate on
$J^{k_\out}(n,1)$.
Consequently, the norm of the second wedge factor in the
r.h.s.\ of \eqref{MinimalNontrivial} is strictly greater than
$|\square_{k_\out}|$, and thence
$\square_{k_\out}[\square_{k_\nel}]$ is trivial.
This completes the proof.
\end{proof}

%
\begin{rem}
The parity of the number $N=\binom{n+k}{n}$ of arguments in
\eqref{NaryVolumeBracket} is arbitrary and hence the reasonings of
Theorem \ref{nDimBaseJacobi} exceed case \eqref{AssocBracket} of
the associative algebras; in particular, for $n=1$ we get Proposition
\ref{J2Zero}, as we claimed in Remark \ref{ForwardProofRem}.
\end{rem}

We give an example of the homotopy $3$-Lie
algebra of polynomials in two variables:

\begin{example}
The space of polynomials
$\mathrm{span}_\Bbbk\langle1,x,y,xy\rangle\subset
\Bbbk_2[x,y]$ endowed with the ternary bracket $\mathbf{1}\wedge
D_x\wedge D_y$ acquires a homotopy $3$-Lie algebra structure.
The commutation relations in this algebra are
$$
[1,x,y]=1,\qquad[1,x,xy]=x,\qquad[1,y,xy]=-y,\quad\text{and}\quad
[x,y,xy]=-xy,
$$
and we see that the structural constants are such that the generators
$x$ and $y$ are mixed.
\end{example}

In this section, we have realized the continualization scheme:
the $N$-ary bracket $W^{0,1,\ldots,N-1}$ is defined on the sequence of
$\Bbbk$-algebras
$$
\Bbbk_N[x]\hookrightarrow\Bbbk[[x]]\hookrightarrow\left[
\genfrac{}{}{0pt}{0}{\Bbbk[[x^1,\ldots,x^n]]}
 { \{\sum_\alpha c_\alpha\cdot x^\alpha\mid\alpha\in\Bbbk,\;
  c_\alpha\in\Bbbk\}. }\right.
$$
The sets of indexes are finite, cardinal, cardinal w.r.t.\ any of $n$
generators, and continuous, respectively, for any possible $N$.

We also note that the definition of the Koszul $\dd_\Delta$-cohomologies
is invariant w.r.t.\ the number of derivations $\dd_i\colon\cA\to\cA$,
$i=1,\ldots,n$, so that the cohomological constructions are preserved
for the Wronskian determinant in \eqref{NaryVolumeBracket}
if $n\geqslant1$.
These concepts allow further, purely algebraic
studies on the topic.
Also, we note that there is a famous mechanism that provides the
associative algebra structures, namely, the Yang\/-\/Baxter
equation (\cite{Kassel}) and the WDVV equation (see~\cite{Geurts} and
references therein).


\section*{Acknowledgements}
The author thanks V.~M.~Buchstaber, I.~S.~Krasil'shchik,
A.~V.~Ovchinnikov, and A.~M.~Verbovetsky
for discussion and remarks, and also A.~A.~Belavin,
B.~L.~Feigin, V.~G.~Kac,
and V.~V.~Sokolov for their attention and advice.
The work was
partially supported by the INTAS grant YS~2001/2-33.

\end{document}